\documentclass [a4paper,twoside,12pt]{article}

\usepackage[utf8]{inputenc}
\usepackage{amsmath,amsfonts,amssymb}
\usepackage{vmargin,graphicx,theorem}
\usepackage[english]{babel}
\usepackage{enumerate}

\setpapersize[portrait]{A4}
\setmarginsrb{2.2cm}{2cm}{2.2cm}{1cm}{0cm}{0.1cm}{0.5cm}{2cm}

\newcommand{\disp}{\displaystyle}

\newcommand{\dN}{\ensuremath{\mathbb{N}}}

\newcommand{\dR}{\ensuremath{\mathbb{R}}}

\newtheorem{ethm}{Theorem}[section]

\newtheorem{ecor}[ethm]{Corollary}

\newtheorem{eprop}[ethm]{Proposition}

\newtheorem{elem}[ethm]{Lemma}

\newtheorem{ex}{Example}

\newtheorem{edefi}[ethm]{Definition}

\theoremstyle{definition}
\newtheorem{remark}{Remark}

\newcommand{\proofend}{~$\rhd$}
\newcommand{\proofbegin}{~$\lhd$}

\newenvironment{eproof}
               {\noindent {\emph{\textbf{Proof}}}\\\proofbegin~}
               {\proofend\\}

\newcommand{\p}[4]{{#3}\!\left#1{#4}\right#2}

\newcommand{\ABS}[1]{\ensuremath{{\left| #1 \right|}}} 
\newcommand{\PAR}[1]{\ensuremath{{\left(#1\right)}}} 
\newcommand{\BRA}[1]{\ensuremath{{\left\{#1\right\}}}} 
\newcommand{\NRM}[1]{\ensuremath{{\left\Vert #1\right\Vert}}} 

\renewcommand{\phi}{\varphi}

\renewcommand{\geq}{\geqslant}

\newcommand{\varf}[1]{\mathbf{Var}_{#1}}

\newcommand{\var}[2]{\p(){\varf{#1}}{#2}}

\def\disp{\displaystyle}

\newcommand{\A}{\ensuremath{\mathcal A}}
\newcommand{\K}{\ensuremath{\mathcal K}}
\newcommand{\dd}{\ensuremath{\mathcal D}}

\newcommand{\R}{\dR}

\newcommand{\C}[1]{\ensuremath{{\mathcal C}^{#1}}}

\newcommand{\PT}[1]{\mathbf{P}_{\!#1}}
\newcommand{\Pt}[1][t]{\ensuremath{\mathbf{P_{\!#1}}}}

\begin{document}

\title{Asymptotic behaviour of  reversible chemical reaction-diffusion equations}
\author{ Ivan Gentil and Boguslaw Zegarlinski }

\date{\today}
\maketitle\thispagestyle{empty}
\abstract{We investigate the asymptotic behavior of the  a large class of  reversible chemical reaction-diffusion equations with the same diffusion. In particular we  prove  the optimal rate in two cases : when there is no diffusion and in the classical "two-by-two" case.
}

\bigskip

\noindent
{\bf Key words:} Reaction-diffusion equation, Spectral gap inequality, Poincar\'e inequality, Exponential decay.  

\bigskip

\noindent
{\bf AMS subject classification:} 35B40, 35A05, 35F25
\bigskip

 \section{Introduction} 
The self-ionization of water  is the chemical reaction in which two water molecules react to produce a hydronium ion $H_3O^+$ and a hydroxide ion $OH^-$. It can be written as follow 
\begin{equation}
\label{eq-ion}
2H_2O\rightleftharpoons H_3O^++OH^-.
\end{equation}
This is a classical  example of a reversible reaction-diffusion chemical reaction between three species such that the water has to appear twice for the reaction.  
\medskip

We would like to investigate in this paper how fast reversible chemical reactions (as the self-ionization) tend to the equilibrium. To be more realistic a diffusion term will be added. Such term models the fact that particles have to move (to diffuse) before the reaction. If the reaction is not enough mixed, then the diffusion term will be more important in the equation. 

\medskip

The main model studied involves $q>1$ species $\A_i$ interacting together as follows 
\begin{equation}
\label{eq-dia}
\sum_{i=1}^q \alpha_i\A_i \rightleftharpoons \sum_{i=1}^q \beta_i\A_i,
\end{equation}
where $\alpha_i,\beta_i\in\dN$, with $\dN$ denoting the set of nonnegative integers. We assume that for any $1\leq i\leq q$, $\alpha_i-\beta_i\neq 0$ which correspond to the case of a reaction without catalyzer, as in the chemical equation presented in~\eqref{eq-ion}.

\medskip

First we will introduce the differential equations that describe the evolution of the species $\A_i$ for $1\leq i\leq q$ in the relation~\eqref{eq-dia}. For some $0\leq i\leq q$ we will denote by $a_i$ the concentration of the species $\A_i$. We assume existence of two non-negative real-valued rate functions 
$$
\K_{\sum_{i=1}^q \alpha_i\A_i \rightarrow \sum_{i=1}^q \beta_i\A_i}, \,\,\K_{\sum_{i=1}^q \beta_i\A_i \rightarrow \sum_{i=1}^q \alpha_i\A_i},
$$
which describe the evolution of the concentration $a_i$ in the two (reversible) reactions.  In the first one we lose $\alpha_i$ molecules of the specie $\A_i$ and in the second one there is a gain of $\beta_i$ molecules of the same species  $\A_i$. We get the opposite  for the reverse reaction.   Thus we can write 
\begin{multline}
\label{eq-ssa}
   \quad\frac{d}{dt}a_i=
   (\beta_i-\alpha_i)\K_{\sum_{i=1}^q \alpha_i\A_i \rightarrow \sum_{i=1}^q \beta_i\A_i}(a)\\-(\beta_i-\alpha_i)\K_{\sum_{i=1}^q \beta_i\A_i \rightarrow \sum_{i=1}^q \alpha_i\A_i}(a), 
\end{multline}
where $a=(a_1,\cdots,a_q)$. We assume the kinetics to be of {\it mass action type}, which means that 
$$
\K_{\sum_{i=1}^q \alpha_i\A_i \rightarrow \sum_{i=1}^q \beta_i\A_i}(a_i)=l\prod_{j=1}^q a_j^{\alpha_j},
$$
where $l$ is a positive constant called {\it the rate constant} of the reaction. This model was proposed by Waage and Guldberg in 1864.  The mass action represents the probability of the reaction between all the species $(\alpha_1,\A_i)_{1\leq i\leq q}$ which are assumed to be independent. Let us denote by $k>0$ the rate constant for the reverse reaction. Of course $l$ and $k$ could be different. 

We will presume that the pot of the reaction is not mixed or not enough mixed. Then  concentrations of species depend on the position in the pot and we have to add a diffusion term which depends also on the species. We obtain that the  model has the following mathematical representation
$$
\partial_t a_i= L_i a_i +(\beta_i-\alpha_i)\PAR{l\prod_{j=1}^q a_j^{\alpha_j} -k\prod_{j=1}^q a_j^{\beta_j}},
$$
where for all $i$, $L_i$ is a diffusion operator. 

\medskip

We refer to \cite{erdi-toth,fife} for a general introduction on chemical reaction-diffusion models. 
\subsection{Mathematical model}
Throughout the entire paper we will assume that all species diffuse with the same speed, i.e. for all $i\in\BRA{1,\cdots,q}$, $L_i=L$, where $L$ is general diffusion generator.  
(While this is a restrictive case - as generally this hypothesis is not realistic - assuming identical diffusion will allow us  to obtain some optimal bounds.)  

In more detail we let $\Omega\subset\dR^n$ ($n\geq1$) to be a bounded open and connected set and assume that the boundary $\partial\Omega$  of $\Omega$,  is $\C{\infty}$-smooth. Let a diffusion operator $L$ be given by
\begin{equation}
\label{eq-defgene}
L f(x)= \sum_{i,j=1}^n {\mathrm a}_{i,j}(x)\partial_{i,j}^2 f(x)+\sum_{i=1}^n b_i(x)\partial_i f(x),
\end{equation}
for any smooth functions $f$, with ${\mathrm a}_{i,j}$ and $b=(b_i)_{1\leq i\leq n}$ in $\C{\infty}(\Omega)$ and the matrix ${\mathbf a}(x)=({\mathrm a}_{i,j}(x))_{i,j}$ symmetric and positive for all $x\in\Omega$. (Generally the choice of coefficients ${\mathrm a}_{i,j}$ and $b_i$ may depend on the domain $\Omega$.) 

We will assume that there exists a probability measure $\mu$ on $\Omega$, $\int\!\! d\mu=1$, which is {\it $L$ - invariant},  i.e. for all functions  $f\in L^1(d\mu)$ in the domain of $L$   
$$
\int Lf\,d\mu=0,
$$
The existence of an invariant measure is a quite standard problem and has a solution for a large class of generators $L$. 

\bigskip

If we denote by $a_i(t,x)$, ($t\geq 0$, $x\in\Omega$), the concentrations of the species $\A_i$ at  time $t$ in the position $x$, then the following reaction-diffusion system is satisfied:
\begin{equation}
\label{eq-def}
  \quad\left\{
\begin{array}{rl}
\disp\forall t> 0,\,\,\forall x\in\Omega,\,\,& \partial_t a_i(t,x) =\disp  L a_i +(\beta_i-\alpha_i)\PAR{l\prod_{j=1}^q a_j^{\alpha_j} -k\prod_{j=1}^q a_j^{\beta_j}},\\
\forall x\in\Omega,\,\,& a_i(0,x)=\disp a_i^0(x)\\
\disp\forall x\in\partial\Omega,\,\, \forall t> 0,\,\,& \disp\frac{\partial a_i }{\partial \nu} (t,x)=0.
\end{array}
\right.
\end{equation}
where $k,l>0$ are rate  constants of the reverse reaction. It is assumed that the initial conditions satisfy  for all $1\leq i\leq q$, $a_i^0\geq 0$ and  $\int a_i^0d\mu>0$. 

The last equation in~(\ref{eq-def}) represents the  Neumann boundary conditions which is quite natural in the context of chemical reaction-diffusion. 

Let us denote by $\dd(\Omega)$ the set of smooth functions $f$ on $\Omega$ satisfying the Neumann boundary conditions, such that  for all $x\in\partial \Omega$, $\frac{\partial f}{\partial \nu} (t,x)=0$.\\   


Let  $(\Pt)_{t\geq}$ be the semi-group associated to $L$, ( which is a linear operator defined  for all functions $f\in\dd(\Omega)$ ), as follows
\begin{equation*}
\left\{
\begin{array}{rl}
\disp\disp\forall t> 0,\,\,\forall x\in\Omega,\,\,&\frac{\partial}{\partial t}\Pt(f)(x)=\disp L\Pt(f)(x)\\
\disp\forall x\in\Omega,\,&\PT{0}(f)(x)=\disp f(x)\\
\disp\forall x\in\partial\Omega,\,\, \forall t> 0,\,\,& \disp\frac{\partial\Pt(f)}{\partial\nu}(x)=0.
\end{array}
\right.
\end{equation*}
Since $\mu$ is an invariant measure, so for all $t\geq0$, $\int \Pt fd\mu=\int fd\mu$.  It could be the classical heat equation with boundary condition but also the heat equation with a drift part some example will be given in the example~\ref{ex-11}.

\medskip

One of the main tools of the article is the spectral gap inequality (also called Poincar\'e inequality). We assume that the operator $L$ satisfies a  {\it spectral gap} inequality in $\dd(\Omega)$, that is : 
There exists  $C_{SG}\geq0$ such that for all smooth functions $f\in\dd(\Omega)$, 
\begin{equation}
\label{eq-poi}
\var{\mu}{f}:=\int\PAR{f-\int f d\mu}^2d\mu \leq -2C_{SG}\int fLf d\mu, 
\end{equation}
which is equivalent to,  for all functions $f\in \dd(\Omega)$ and $t\geq 0$:
\begin{equation}
\label{eq-sg}
\var{\mu}{\Pt f}\leq e^{-\frac{1}{C_{SG}}t}\var{\mu}{f}.
\end{equation}
In the general form, when $L$ is given by ~\eqref{eq-defgene}, then for all function $f$ with Neumann boundary conditions, 
$$
-\int fLfd\mu = \int \nabla f (x) \cdot({\mathbf a}(x) \nabla f(x))d\mu(x),
$$
where ${\mathbf a}(x)$ is the diffusion matrix of the generator $L$ given in~\eqref{eq-defgene}.
The assumption of the spectral gap inequality is fundamental. First able it proves that the semigroup $(\Pt)_{t\geq0}$ is {\it ergodic}, that is for all $f\in\mathcal D(\Omega)$,
\begin{equation}
\label{eq-ergodi}
\lim_{t\rightarrow\infty}\Pt f=\int f d\mu,
\end{equation}
in $L^2(d\mu)$. Moreover this inequality proves that  the rate of convergence is exponential in $L^2(d\mu)$, it is given by the inequality~\eqref{eq-sg}. See for instance chapter~2 of  \cite{logsob} of a review on semigroups tools and ergodic properties. 
\begin{ex}
\label{ex-11}
Here are some basic examples on a bounded domain $\Omega$.

One can consider the generator $L$ given in~\eqref{eq-defgene} where  $b=0$ and  ${\mathbf a}\geq\lambda Id$ in $\Omega$ with $\lambda>0$ in the sense of symmetric matrices.  If ${\mathbf a}=Id$ then the associated semigroup is $(\Pt)_{t\geq0}$ is the classical heat semigroup with Neumann boundary conditions. The invariant measure is the Lebesgue measure and the spectral gap constant $C_{SG}$, depends on the domain~$\Omega$.  

But we can also consider a generator $L=\Delta -\nabla \psi\nabla $ with $\psi$ a smooth function on the domain $\Omega$. The generator  $L$ is symmetric in $L^2(\mu_\psi)$ where $\mu_\psi=e^{-\psi}dx/C$  is a probability measure in $\Omega$ with the normalization constant $C$.  It satisfies a Spectral gap  inequality with the dirichlet form $-\int\!  fLfd\mu_\psi=\int\! |\nabla f|^2d\mu_\psi $ for  functions $f$  satisfying the  Neumann Boundary conditions. 
\end{ex}

\medskip

We will first study  the problem of a steady state of the differential system~(\ref{eq-def}).  Let $v_i(t,x)=\lambda_i a_i( t, x)$ for some constant $\lambda_i>0$, assuming that $a_i$ is solution of~(\ref{eq-def}), then $v_i$ satisfies 
$$
\partial_t v_i=L v_i +{ \lambda_i}(\beta_i-\alpha_i)\PAR{l\prod_{i=1}^q \PAR{\frac{ v_i}{\lambda_i}}^{\alpha_i} -k\prod_{i=1}^q \PAR{\frac{ v_i}{\lambda_i}}^{\beta_i}}. 
$$
Let us chose constants $\lambda_i>0$ such that $k\prod_{i=1}^q (\lambda_i)^{\alpha_i} =l\prod_{i=1}^q (\lambda_i)^{\beta_i}$, or equivalently  $\prod_{i=1}^q (\lambda_i)^{\alpha_i-\beta_i} =k/l$. We then obtain that $v_i$ is solution of 
\begin{equation}
\label{eq-def2}
  \,\left\{
\begin{array}{rl}
\disp\forall t> 0,\,\,\forall x\in\Omega,\,\,&\partial_t v_i(t,x)= L v_i(t,x) +k_i(\beta_i-\alpha_i)G(v_1(t,x),\cdots,v_q(t,x)),\\
\forall x\in\Omega,\,\,& v_i(0,x)=\disp v_i^0(x)\\
\disp\forall x\in\partial\Omega,\,\, \forall t> 0,\,\,&\disp\frac{\partial v_i}{\partial\nu}(t,x)=0
\end{array}
\right.
\end{equation}
where $G(v_1,\cdots,v_q)={\prod_{j=1}^q v_j^{\alpha_j} -\prod_{j=1}^q v_j^{\beta_j}}$ and 
\begin{equation}
\label{eq-ki}
k_i=\frac{\lambda_i l}{\prod_{i=1}^q (\lambda_i)^{\alpha_i}}. 
\end{equation}
%

\noindent Let us define the set
$\mathcal{S}=\displaystyle\BRA{(z_i)_{1\leq i\leq q},\,\,{\rm{s.t.}}\,\,\sum_{i=1}^q z_i k_i(\beta_i-\alpha_i)=0 }$,
 where $k_i$ is defined in~(\ref{eq-ki}). Then for all $(z_i)_{1\leq i\leq q}\in \mathcal{S}$, one gets 
$$
\partial _t \sum_{i=1}^q z_i v_i=L\sum_{i=1}^q z_i v_i,
$$
which gives 
$$
\sum_{i=1}^q z_i v_i(t,x)=\Pt\PAR{\sum_{i=1}^q z_i v_i^0}(x). 
$$
In particular, due to the fact that $\mu$ is an invariant measure, $$
\sum_{i=1}^q z_i\int v_i(t,x)d\mu(x)=\sum_{i=1}^q z_i\int v_i^0d\mu.
$$

The goal of this article is to understand the asymptotic behaviour of the reaction-diffusion problem in the spirit of \cite{df06,df08}, for a general but (diagonal) diffusion. 

We will first define  a {\it Steady State} of equation~(\ref{eq-def2}), in the following way  standard  in a chemical reversible reaction.  
\begin{edefi}
A {\it steady state} of equation~(\ref{eq-def2}) with non-negative initial conditions $\PAR{v_i^0}_{1\leq i\leq q}$ is   a vector  $(s_i)_{1\leq  i\leq q}\in\PAR{\dR^{+}}^q$  such that for all $(z_i)_{1\leq i\leq q}\in\mathcal{S}$: 
$$
\sum_{i=1}^q z_is_i=\sum_{i=1}^q z_i\int v_i^0d\mu\,\,\,\,
\rm{and} \,\,\,\,
{\prod_{i=1}^q s_i^{\alpha_i} =\prod_{i=1}^q s_i^{\beta_i}}.
$$
\end{edefi}

\begin{remark}
We implicitly  assume in the previous definition that a steady state is a vector  $(s_i)_{1\leq  i\leq q}\in\PAR{\dR^{+}}^q$ independent of $x$. In fact we don't know in the general case if there exists a steady state depending of the space variable $x$.  But as it is proved in this paper, the solution converges to the one defined above. 
\end{remark}
\begin{elem}
\label{lem-sta}
Let $(v_i)_{1\leq i\leq q}$ satisfies equation~$(\ref{eq-def2})$ with initial conditions satisfying $v_i^0\geq 0$ and  $\int v_i^0d\mu> 0$. Then there exists a unique steady  state $(s_i)_{1\leq  i\leq q}$ of~$(\ref{eq-def2})$  such that for all $i\in\BRA{1,\cdots,q}$,  $s_i>0$.
\end{elem}

\begin{eproof}
The steady state $(s_i)_{1\leq i\leq n}$ has to satisfy:
$$
\forall (z_i)_{1\leq i\leq q}\in \mathcal{S},\,\, \sum_{i=1}^q z_is_i=\sum_{i=1}^q z_i\int v_i^0d\mu:=M_z. 
$$
Let $A=\BRA{(s_i)_{1\leq i\leq q},\,\, \sum_{i=1}^q z_is_i=M_z}$. Since for all $i$, $\alpha_i\neq\beta_i$ then  $\mathcal{S}$ is a subspace of $\dR^q$ of dimension $q-1$,  then $A$ is a manifold of dimension 1. Then one gets
$$
A=\BRA{\PAR{\int v_i^0d\mu+t k_i(\beta_i-\alpha_i)}_{1\leq i\leq q},\,\,t\in\dR}.
$$
We have to find $t\in\dR$, such that for all  $1\leq i\leq q$, $\int v_i^0d\mu+tk_i(\beta_i-\alpha_i)>0$  and 
\begin{equation*}
\label{eq-st}
\prod_{i=1}^q \PAR{\int v_i^0d\mu+tk_i(\beta_i-\alpha_i)}^{\beta_i-\alpha_i} =1.
\end{equation*}
A simple computation gives that the function 
$$
\phi(t)=\prod_{i=1}^q\PAR{\int v_i^0d\mu+tk_i(\beta_i-\alpha_i)}^{\beta_i-\alpha_i}
$$
is defined on the set $[a,b)$ with 
\begin{equation*}
\left\{
\begin{array}{r}
\displaystyle a=\inf\BRA{t,\,\,{\rm{s.t.}}\,\,\forall i,\,\int v_i^0d\mu+tk_i(\beta_i-\alpha_i)\geq0}\\
\displaystyle b=\sup\BRA{t,\,\,{\rm{s.t.}}\,\,\forall i,\,\int v_i^0d\mu+tk_i(\beta_i-\alpha_i)\geq0}. 
\end{array}
\right.
\end{equation*}
For all $t\in[a,b)$, one gets 
$$
\frac{\phi'(t)}{\phi(t)}=\sum_{i=1}^q\frac{k_i(\beta_i-\alpha_i)^2}{\int v_i^0d\mu+tk_i(\beta_i-\alpha_i)},
$$
then $\phi$ is increasing and satisfies $\phi(a)=0$ and $\phi(b)=+\infty$. Thus there exists a unique $t\in(a,b)$ such that $\phi(t)=1$, which ends the proof. 
\end{eproof}


\medskip

Due to the fact that all species are moving according to the same diffusion, one can reduce the problem as follow: for all $1\leq i,j\leq q$, one gets
$$
\partial_t\PAR{\frac{v_i}{k_i(\beta_i-\alpha_i)}-\frac{v_j}{k_j(\beta_j-\alpha_j)}}=L\PAR{\frac{v_i}{k_i(\beta_i-\alpha_i)}-\frac{v_j}{k_j(\beta_j-\alpha_j)}}
$$
and thus 
\begin{equation}
\label{eq-ct}
\begin{array}{rl}
\displaystyle v_i& =\displaystyle\frac{k_i(\beta_i-\alpha_i)}{k_j(\beta_j-\alpha_j)}v_j+\Pt\PAR{v_i^0-\frac{k_i(\beta_i-\alpha_i)}{k_j(\beta_j-\alpha_j)}v_j^0}\\
&:=\displaystyle\frac{k_i(\beta_i-\alpha_i)}{k_j(\beta_j-\alpha_j)}v_j+C_{i,j}(t,x). 
\end{array}
\end{equation}

Let us fix $1\leq i\leq q$. Then the study of~(\ref{eq-def2}) is equivalent  to the study of the following PDE with boundary conditions
\begin{equation}
\label{eq-apt}
  \begin{array}{rl}
\disp\partial_t v_i& =\disp  L v_i +k_i(\beta_i-\alpha_i)\times\\
&\!\!\!\!\disp\PAR{\prod_{j=1}^q \PAR{\frac{k_j(\beta_j-\alpha_j)}{k_i(\beta_i-\alpha_i)}v_j+C_{j,i}(t,x)}^{\alpha_j} -\prod_{j=1}^q \PAR{\frac{k_j(\beta_j-\alpha_j)}{k_i(\beta_i-\alpha_i)}v_j+C_{j,i}(t,x)}^{\beta_j}}\\
&:=\disp L v_i +k_i(\beta_i-\alpha_i)F(t,x,v_i).
\end{array} 
\end{equation}

The existence problem in the general case when the operators $L_i$ depending on~$i$ are all different is a difficult problem. Some approaches one can find e.g. in \cite{feng,dfpv,gv,pierre,smoller} and also by discretization in~\cite{glitzky,glitzky2}. 

\medskip

Our contribution to the domain of reaction-diffusion equation is to prove, using  the spectral gap inequality (or Poincar\'e inequality) and Markov semigroup tools,  that the solution of reaction-diffusion equation~(\ref{eq-def2})  converges to the unique steady state associated to the initial condition. The goal of this paper is to explain how fast  the solution converges to the steady state and to debate on  the speed of convergence. In particular if the result obtained is far or close  to the optimal one. 

The main idea is to investigate a simple proof to explain the asymptotic behavior. Therefore, for simplicity, we will not focus on existence theorems, the optimal conditions of initial conditions and do not study the special case of a unbounded space $\Omega$. Almost all results given here can be generalize to the entire space $\dR^n$ if a regular solution of the problem is given.

\medskip

We will consider different cases. In Section~\ref{sec-11} we will study the case  without diffusion,  that means that the concentration of different species do not depend on $x\in\Omega$. We will prove that there exists a solution converging to the unique steady state with an exponential and optimal rate of convergence.

Then in Section~\ref{sec-22} we study the classical  ``two-by-two'' case  
$$
\A+\mathcal B\rightleftharpoons \mathcal C +\mathcal D,
$$
treated in \cite{df08,dfm} by entropy methods. This case is interesting because we obtain the optimal rate of convergence and  it gives  tools to understand a more general case.

We finish in Section~\ref{sec-5} with the general case
$$
  \sum_{i=1}^q \alpha_i\A_i \rightleftharpoons \sum_{i=1}^q \beta_i\A_i.
$$
We prove that under the existence of a non-negative solution and under the assumption that $\alpha_i\beta_i=0$, we get that solution converges with an exponential rate to the steady state. In this case the rate may be not optimal.

 \section{Case without diffusion}
 \label{sec-11}
Assume that concentrations of species do not depend on $x\in\Omega$. It is the case the pot used for the chemical reaction is mixed constantly so that its contents remain spatially homogeneous.  This case is important because we can solve it explicitly and it gives tools to study the general case in the later sections. 

We believe  that results of this section are not new, we give here the proof to keep the article self contained and also because we did not find appropriate reference. 

\medskip

The chemical reaction, without diffusion is given by the following system, for all $i\in\BRA{1,\cdots,q}$,
\begin{equation}
\label{eq-rea}
\frac{d}{dt} v_i=k_i(\beta_i-\alpha_i)\PAR{\prod_{j=1}^q v_j^{\alpha_j} -\prod_{j=1}^q v_j^{\beta_j}}, 
\end{equation}
where $k_i$ is defined in~(\ref{eq-ki}), with initial conditions  $v_i(0)>0$ for $1\leq i\leq q$. Equivalently, using the same method as for~(\ref{eq-apt}),  for some  $1\leq i\leq q$ fixed, 
\begin{multline*}
\frac{d}{dt} v_i= \\\disp k_i(\beta_i-\alpha_i)\PAR{\prod_{j=1}^q \PAR{\frac{k_j(\beta_j-\alpha_j)}{k_i(\beta_i-\alpha_i)}v_i+C_{j,i}}^{\alpha_j} -\prod_{j=1}^q \PAR{\frac{k_j(\beta_j-\alpha_j)}{k_i(\beta_i-\alpha_i)}v_i+C_{j,i}}^{\beta_j}}\\
:=\disp k_i(\beta_i-\alpha_i)F(v_i),
\end{multline*}
where  
\begin{equation}
\label{eq-defcji}
C_{j,i}=k_j\PAR{\beta_j-\alpha_j}\PAR{\frac{v_j(0)}{k_j(\beta_j-\alpha_j)}-\frac{v_i(0)}{k_i(\beta_i-\alpha_i)}}.
\end{equation}

\begin{ethm}
\label{thm-reaction}
Let the initial conditions $(v_j(0))_{1\leq j\leq q}$ be positive.  Then equation~$(\ref{eq-rea})$ has a unique solution defined on $[0,+\infty)$, which satisfies for all $1\leq i\leq q$,
\begin{equation*}
\ABS{v_i(t)-s_i}\leq e^K\ABS{v_i(0)-s_i}\exp\PAR{-Ct},
\end{equation*} 
where $K$ is a constant depending on initial conditions, the steady state $(s_i)_{1\leq i\leq q}$ is defined in Lemma~\ref{lem-sta} and
\begin{equation}
\label{eq-defc}
C=\prod_{i=1}^q {s_i}^{\alpha_i}\sum_{i=1}^q \frac{k_i(\beta_i-\alpha_i)^2}{s_i}.
\end{equation}
Moreover the constant  $C$ is the optimal rate of convergence.  
\end{ethm}

\begin{eproof}
Let $j_0$ be such that $\beta_{j_0}-\alpha_{j_0}>0$ and for all $i$ s.t. $\beta_{i}-\alpha_{i}>0$, one has $\frac{v_{j_0}(0)}{k_{j_0}(\beta_{j_0}-\alpha_{j_0})}\leq \frac{v_{i}(0)}{k_i(\beta_{i}-\alpha_{i})}$. If the set $\BRA{i;\,\,\beta_{i}-\alpha_{i}>0 }$ is empty one can use the negative part. Assume for simplicity that $j_0=1$. Then the reaction equation becomes
\begin{multline}
\label{eq-sim}
  \frac{d}{dt} v_1= k_1(\beta_1-\alpha_1)\times\\ 
\!\!\!\!\PAR{\prod_{i=1}^q \PAR{\frac{k_i(\beta_i-\alpha_i)}{k_1(\beta_1-\alpha_1)}v_1+C_{i,1}}^{\alpha_i} -\prod_{i=1}^q \PAR{\frac{k_i(\beta_i-\alpha_i)}{k_1(\beta_1-\alpha_1)}v_1+C_{i,1}}^{\beta_i}}.
\end{multline}
By the definition of $j_0$ one have, $C_{i,1}\geq0$ if $\beta_i-\alpha_i>0$ and $C_{i,1}\leq0$ if $\beta_i-\alpha_i<0$.

Getting a positive solution $(v_i)_{1\leq i\leq q)}$ of~(\ref{eq-rea}) is equivalent to getting a solution $v_1$ of~(\ref{eq-sim}) which is defined by 
\begin{equation*}
\displaystyle v_i =\displaystyle\frac{k_i(\beta_i-\alpha_i)}{k_j(\beta_j-\alpha_j)}v_j+{v_i^0-\frac{k_i(\beta_i-\alpha_i)}{k_j(\beta_j-\alpha_j)}v_j^0}
\end{equation*}
then it satisfies the following inequality
$$
\forall t\in[0,+\infty),\quad 0<v_1(t)<M,
$$
where $M=v_1(0)-k_1(\beta_1-\alpha_1)\max\BRA{\frac{v_i(0)}{k_i(\beta_i-\alpha_i)};\,\,\beta_i-\alpha_i<0}$. By convention, if the set  $\BRA{i;\,\,\beta_i-\alpha_i<0}$ is empty,  then we have $\max\BRA{\frac{v_i(0)}{k_i(\beta_i-\alpha_i)};\,\,\beta_i-\alpha_i<0}=-\infty$.

Let us denote by
\begin{multline}
\label{eq-deff}
   F(X)=k_1(\beta_1-\alpha_1)\times\\
\!\!\!\!\PAR{\prod_{i=1}^q \PAR{\frac{k_i(\beta_i-\alpha_i)}{k_1(\beta_1-\alpha_1)}X+C_{i,1}}^{\alpha_i} -\prod_{i=1}^q \PAR{\frac{k_i(\beta_i-\alpha_i)}{k_1(\beta_1-\alpha_1)}X+C_{i,1}}^{\beta_i}}.
\end{multline}
Lemma~\ref{lem-sta} proves that the polynomial equation  $F=0$ has only one solution $s_1$ in the set $(0,M)$. Let $Q$ be a factor of $F$, i.e. we have a factorization $F(X)=(X-s_1)Q(X)$. Then a simple computation yields 
\begin{multline*}
   F'(s_1)=-\prod_{i=1}^q \PAR{\frac{k_i(\beta_i-\alpha_i)}{k_1(\beta_1-\alpha_1)}s_1+C_{i,1}}^{\alpha_i}\sum_{i=1}^q \frac{k_i(\beta_i-\alpha_i)^2}{\frac{k_i(\beta_i-\alpha_i)}{k_1(\beta_1-\alpha_1)}s_1+C_{i,1}}=\\
\quad\quad\quad\quad\quad\quad\quad\quad\quad\quad\quad\quad\quad-\prod_{i=1}^q {s_i}^{\alpha_i}\sum_{i=1}^q \frac{k_i(\beta_i-\alpha_i)^2}{s_i}:=-C<0,
\end{multline*}
which proves that $F'(s_1)=Q(s_1)<0$ and then $Q(X)<0$ for all $X\in(0,M)$.

Function $F$ is locally a Lipschitz function, thus by Cauchy-Lipschitz's Theorem, there exists a unique maximal solution starting at $v_1(0)$ of the equation 
$$
\forall t\in[0,T),\,\,\frac{dv_1(t)}{dt}=F(v_1(t))=(v_1(t)-s_1)Q(v_1(t)),
$$ 
for some $T>0$. 

One has $Q(X)<0$ for all $X\in(0,M)$, which implies that if $v_1(0)>s_1 $, then $v_1$ is non-increasing and moreover $v_1(t)\geq s_1$ for all $t\in[0,T[$, while if $v_1(0)<s_1 $, then $v_1$ is non-decreasing and moreover $v_1(t)\leq s_1$ for all $t\in[0,T[$, and if $v_1(0)=s_1$, then $v_1(t)=s_1$ for all $t\in[0,T[$. Then one gets that $T=+\infty$ and for all $t\geq0$, $v_1(t)\in(0,M)$. 

Using the identity
$$
\frac{1}{F(X)}=\frac{1/Q(s_1)}{X-s_1}+\frac{R(X)}{Q(X)},
$$
with $R(X)=(Q(s_1)-Q(X))/(Q(s_1)(X-s_1))$ we get  for all $t\geq 0$
\begin{equation}
\label{eq-so}
  \quad\quad\quad\quad\ABS{v_1(t)-s_1}=\ABS{v_1(0)-s_1}\exp\PAR{Q(s_1)t+\int_{v_1(0)}^{v_1(t)}\frac{Q(a)-Q(s_1)}{(a-s_1)Q(a)}da}.
\end{equation}
This identity  proves that $v_1$ goes to $s_1$ and moreover if 
$$
K=\int_{v_1(0)}^{s_1}\frac{Q(a)-Q(s_1)}{(a-s_1)Q(a)}da,
$$
then for all $t\geq0$,
$$
\ABS{v_1(t)-s_1}\leq e^{|K|}\ABS{v_1(0)-s_1}\exp\PAR{-Ct}.
$$
Links between $v_i$ and $v_1$ gives the last inequalities for any $v_i$.  Moreover, equality~(\ref{eq-so}) proves that the constant $C$ is the optimal rate. 
\end{eproof}

 \section{The ``two-by-two'' case}
 \label{sec-22}
The goal of this section is to investigate the asymptotic behaviour of a  chemical reaction of particular type
$$
\A+\mathcal B\rightleftharpoons \mathcal C +\mathcal D.
$$
We will assume that all species are moving according to the same generator and the speed of the two reactions are the same, for instance equal to 1.

\medskip

This case was treated with a general diffusion in~\cite{df08,dfm} by using entropy method, but the optimal rate of convergence was not obtained. 

\medskip

 Let us denote by $a$, $b$, $c$ and $d$ concentrations of $\mathcal A$,   $\mathcal B$, $\mathcal C$ and  $\mathcal D$ in the domain $\Omega$. In this case the functions $a$, $b$, $c$ and $d$ are solutions of the following system on $\Omega$,
\begin{equation}
\label{eq-nlin}
\left\{
\begin{array}{r}
\partial_t a= L a -ab+cd\\
\partial_t b= L b -ab+cd\\
\partial_t c= L c +ab-cd\\
\partial_t d= L d +ab-cd\\
\end{array}
\right.
\end{equation}
with non-negative initial conditions   $a_0$, $b_0,$ $c_0$ and $d_0$ in $\mathcal D(\Omega)$, such that $\int a_0d\mu>0$ and the same for   $b_0,$ $c_0$ and $d_0$. We recall that the boundary conditions are included in the definition of the domain $\mathcal D(\Omega)$. 

If $a,b,c,d$ are solutions, then  $a+c=\Pt\PAR{a_0+c_0}$, $a+d=\Pt\PAR{a_0+d_0}$ and $a-b=\Pt\PAR{a_0-b_0}$ which gives that the function $a$ is solution of the linear equation
\begin{equation}
\label{eq-derder2}
\partial _t a=La-aD_t+C_t,
\end{equation}
where $C_t=\Pt(a_0+c_0)\Pt(a_0+d_0)$ and $D_t=\Pt\PAR{a_0+b_0+c_0+d_0}$.  

Since the solution is given by a linear solution there exists a classical non-negative solution of the problem. One can see for example \cite{amann1,rothe} for a proof . Some remarks on the existence are given in Section~\ref{sec-5}.

\medskip

We also obtain the estimation useful for the asymptotic behaviour. 
\begin{elem}
Solution $a$, $b$, $c$ and $d$ of~$(\ref{eq-nlin})$ satisfy  for all $t\geq0$ and $x\in\Omega$,
\begin{equation}
\label{eq-sys}
\left\{
\begin{array}{l}
0\leq a(t,x)\leq \min\BRA{\Pt(a_0+c_0),\Pt(a_0+d_0)}\\
0\leq b(t,x)\leq \min\BRA{\Pt(b_0+c_0),\Pt(b_0+d_0)}\\
0\leq c(t,x)\leq \min\BRA{\Pt(c_0+a_0),\Pt(c_0+b_0)}\\
0\leq d(t,x)\leq \min\BRA{\Pt(d_0+a_0),\Pt(d_0+b_0)}
\end{array}
\right.
\end{equation}
\end{elem}
%




\begin{ethm}
\label{thm-2}
Assume that the semigroup $(\Pt)_{t\geq 0}$  satisfies a spectral gap inequality~$(\ref{eq-sg})$ with respect to the invariant probability measure $\mu$. Let $a_0$, $b_0$, $c_0$ and $d_0$ be non-negative initial conditions satisfying $a_0,b_0,c_0,d_0\in L^{4}(d\mu)$. 

We set $M_{a+b+c+d}=\int (a_0+b_0+c_0+d_0)d\mu$, $M_4=\PAR{\int\PAR{{a_0+b_0+c_0+d_0}}^{4}d\mu}^{\frac{1}{2}}$ and $C_{SG}$ denotes the constant in the spectral gap inequality~$(\ref{eq-sg})$.  

Let $s_a$ is the steady state associated to the initial condition $a_0$, if $M_{a+b+c+d}\neq\frac{1}{8C_{SG}}$, then
the solution $a$ of~$\eqref{eq-derder2}$ satisfies, for all $t\geq0$, 
\begin{multline}
\label{eq-thm2}
   \sqrt{\int\PAR{a-s_a}^2 d\mu}\leq
\PAR{\sqrt{\int\PAR{a_0-s_a}^2 d\mu}+ \ABS{\frac{5M_4}{M_{a+b+c+d}-\frac{1}{8C_{SG}}}}}\times\\
\exp\PAR{-\min\BRA{M_{a+b+c+d},\frac{1}{8C_{SG}}}t},
\end{multline}
and  if $M_{a+b+c+d}=\frac{1}{8C_{SG}}$ then for all $t\geq0$, 
$$
\sqrt{\int\PAR{a-s_a}^2 d\mu}\leq 
\PAR{\sqrt{\int\PAR{a_0-s_a}^2 d\mu}+ 5M_4\,t}\exp\PAR{-M_{a+b+c+d}t}. 
$$
The same inequality holds for  $b$, $c$ and $d$ associated to $s_b$, $s_c$ and $s_d$. 

\medskip

If the initial conditions satisfy  $M_{a+b+c+d}<\frac{1}{8C_{SG}}$, then the rate of the convergence is optimal. 
\end{ethm}

Let us start with a general estimate. 
\begin{elem}
\label{lem-3}
Assume that the semigroup $(\Pt)_{t\geq 0}$  satisfies a spectral gap inequality~$(\ref{eq-sg})$ with respect to the invariant probability measure $\mu$, then for all functions $f\in L^4(\mu)$ and  all $t\geq 0$, 
\begin{equation}
\label{eq-po4}
\int \PAR{\Pt f-\int fd\mu}^4d\mu\leq 4e^{-\frac{1}{2C_{SG}}t}\int f^4 d\mu.
\end{equation}
\end{elem}

\begin{eproof}
Set $\tilde{f}=f-\int fd\mu$, then using semigroup properties and the Cauchy-Schwartz inequality applied to $(\Pt)_{t\geq 0}$, one gets,  
$$
\int \PAR{\Pt f-\int fd\mu}^4d\mu\!=\!\int \PAR{\Pt \tilde{f}}^4d\mu=\int \PAR{\PT{\frac{t}{2}}\PT{\frac{t}{2}} \tilde{f}}^4d\mu, 
$$
since $\PT{\frac{t}{2}}\PT{\frac{t}{2}} \tilde{f}=\PT{{t}}\tilde{f}$. Now the Markov semigroup $(\Pt)_{t\geq0}$ is given by a Markov kernel, so one has
$ \PAR{\PT{\frac{t}{2}}\tilde{f}}^2\leq  {\PT{\frac{t}{2}}(\tilde{f}^2)}$, which gives at the end 
$$
\int \PAR{\Pt f-\int fd\mu}^4d\mu\!\leq
\int \PAR{\PT{\frac{t}{2}} \PAR{\PT{\frac{t}{2}}\tilde{f}}^2}^2 d\mu. 
$$
If we set  $F=\PT{\frac{t}{2}}(\tilde f)$, then one has
$$
\int \PAR{\Pt \tilde{f}}^4d\mu\leq 2 \int\PAR{\PT{\frac{t}{2}}(F^2)-\int F^2 d\mu }^2d\mu+2\PAR{\int {F^2}d\mu}^2,
$$
which gives by definition of $F$,
$$
\int \PAR{\Pt \tilde{f}}^4d\mu\leq 2 \int\PAR{\PT{\frac{t}{2}}(F^2)-\int F^2 d\mu }^2d\mu+2\PAR{\int \PAR{\PT{\frac{t}{2}}(f)-\int fd\mu}^2d\mu}^2.
$$
We apply twice inequality (\ref{eq-sg}) to $F$ and to $f$ to obtain
$$
\int \PAR{\Pt \tilde{f}}^4d\mu\leq 2e^{-\frac{t}{2C_{SG}}}\var{\mu}{F^2}+2e^{-\frac{t}{C_{SG}}}\PAR{\var{\mu}{f^2}}^2,
$$
which implies~(\ref{eq-po4}). 
\end{eproof}

{\noindent {\emph{\textbf{Proof of Theorem~\ref{thm-2}}}}\\\proofbegin~}
Is this  case, the steady state is the following limit,  
\begin{equation}
\label{eq-sa}
s_a=\frac{\int\PAR{a_0+c_0}d\mu\int\PAR{a_0+d_0}d\mu}{\int\PAR{a_0+b_0+d_0+c_0}d\mu}=
\lim_{t\rightarrow +\infty}\frac{C_t}{D_t},
\end{equation}
the limit can be seen in $L^4(d\mu)$. Let us denote by  $M_{a+c}=\int (a_0+c_0)d\mu$ and define similarly  $M_{a+d}$ and $M_{a+b+c+d}$. 
One has, 
$$
\frac{d}{dt}\frac{1}{2}\int\PAR{a-s_a}^2 d\mu=\int (a-s_a)\partial _tad\mu, 
$$
then by~\eqref{eq-derder2} and~\eqref{eq-sa} one obtains 
\begin{eqnarray}
  \frac{d}{dt}\frac{1}{2}\int\PAR{a-s_a}^2 d\mu=\int \!a\, L a\,d\mu-M_{a+b+c+d}\int (a-s_a)^2d\mu\nonumber\\
\!\!\!\!\!\!\!\!\!+\int a(a-s_a)(M_{a+b+c+d}-D_t)d\mu+\int (a-s_a)(C_t-M_{a+c}M_{a+d})d\mu.\label{eq-dt}
\end{eqnarray}
Let us consider the last term:
\begin{eqnarray*}
  {\int (a-s_a)(C_t-M_{a+c}M_{a+d})d\mu}={\int (a-s_a)M_{a+d}(\Pt(a_0+c_0)-M_{a+c})d\mu}\\
\quad\quad\quad\quad\quad\quad\quad\quad+{\int (a-s_a)\Pt(a_0+c_0)(\Pt(a_0+d_0)-M_{a+d})d\mu}.
\end{eqnarray*}
Setting $\phi(t)=\sqrt{\int\PAR{a-s_a}^2 d\mu}$, Cauchy-Schwarz inequality yields
\begin{eqnarray*}
  {\int (a-s_a)(C_t-M_{a+c}M_{a+d})d\mu}\leq\phi(t)\PAR{M_{a+d}\sqrt{\int (\Pt(a_0+c_0)-M_{a+c})^2d\mu}}\\
+\phi(t)\PAR{\PAR{\int \Pt(a_0+c_0)^4d\mu}^{1/4}\PAR{\int (\Pt(a_0+d_0)-M_{a+d})^4d\mu}^{1/4}}.
\end{eqnarray*}
First spectral gap inequality gives 
$$
\int (\Pt(a_0+c_0)-M_{a+c})^2d\mu\leq e^{-\frac{1}{C_{SG}}t}\var{\mu}{a_0+c_0}\leq e^{-\frac{1}{C_{SG}}t}\int(a_0+c_0)^2d\mu.
$$
Since the semigroup $(\Pt)_{t\geq0}$ is contractive : $\frac{d}{dt}\int \Pt(a_0+b_0)^4d\mu\leq 0$, one obtains  
$$
{\PAR{\int \Pt(a_0+c_0)^4d\mu}^{1/4}}\leq{\PAR{\int (a_0+c_0)^4d\mu}^{1/4}}. 
$$
To finish, Lemma~\ref{lem-3} gives
$$
\PAR{\int (\Pt(a_0+d_0)-M_{a+d})^4d\mu}^{1/4}\leq \sqrt{2}e^{-\frac{t}{8C_{SG}}}\PAR{\int (a_0+d_0)^4 d\mu}^{1/4},
$$
which implies for the last term of~(\ref{eq-dt}): 
\begin{multline*}
  {\int (a-s_a)(C_t-M_{a+c}M_{a+d})d\mu}\\
  \leq
3\phi(t)\exp\PAR{-\frac{t}{8C_{SG}}}\PAR{\int\PAR{{a_0+b_0+c_0+d_0}}^{4}d\mu}^{\frac{1}{2}}
\end{multline*}
For the other term one gets
\begin{multline*}
  \int a(a-s_a)(M_{a+b+c+d}-D_t)d\mu\\\leq\phi(t) \PAR{\int a^4 d\mu}^{1/4}\PAR{\int\PAR{M_{a+b+c+d}-D_t}^4d\mu}^{1/4}.
\end{multline*}
Using~(\ref{eq-sys}), we get  
$$
\PAR{\int a^4 d\mu}^{1/4}\leq\min\BRA{\PAR{\int (a_0+c_0) ^4d\mu}^{1/4},\PAR{\int (a_0+d_0) ^4d\mu}^{1/4} },
$$
and~(\ref{eq-po4}) gives
\begin{multline*}
  \int a(a-s_a)(M_{a+b+c+d}-D_t)d\mu\\\leq 2 \phi(t) \PAR{\int\PAR{{a_0+b_0+c_0+d_0}}^{4}d\mu}^{\frac{1}{2}}\exp\PAR{-\frac{t}{8C_{SG}}}.
\end{multline*}
Then we obtain 
$$
\phi'(t)\leq-M_{a+b+c+d}\phi(t)+5M_4\exp\PAR{-\frac{t}{8C_{SG}}t},
$$
where $M_4=\PAR{\int\PAR{{a_0+b_0+c_0+d_0}}^{4}d\mu}^{\frac{1}{2}}$. Integration of the last differential inequality yields:\\
if $M_{a+b+c+d}\neq\frac{1}{8C_{SG}}$, then 
\begin{eqnarray*}
  \sqrt{\int\PAR{a-s_a}^2 d\mu}\leq 
\PAR{\sqrt{\int\PAR{a_0-s_a}^2 d\mu}+ \ABS{\frac{5M_4}{M_{a+b+c+d}-\frac{1}{8C_{SG}}}}}\times\\
\quad\quad\quad\quad\quad\quad\quad\quad\quad\quad\quad\exp\PAR{-\min\BRA{M_{a+b+c+d},\frac{1}{8C_{SG}}}t},
\end{eqnarray*}
and if $M_{a+b+c+d}=\frac{1}{8C_{SG}}$, then 
$$
\sqrt{\int\PAR{a-s_a}^2 d\mu}\leq 
\PAR{\sqrt{\int\PAR{a_0-s_a}^2 d\mu}+ 5M_4\,t}\exp\PAR{-M_{a+b+c+d}t},
$$
which finished the proof of~(\ref{eq-thm2}).

If $M_{a+b+c+d}<\frac{1}{8C_{SG}}$, then the rate becomes $e^{-M_{a+b+c+d}t}$, one can check  that $M_{a+b+c+d}$ is equal to the constant $C$ defined in~(\ref{eq-defc}) which is optimal. 
{\proofend}

\begin{remark}
In the case of a linear equation of diffusion, $\partial_t u=L u$ where $L$ is given by~\eqref{eq-defgene}, the optimal rate of convergence in $L^2(d\mu)$ is given by the spectral gap constant, which is the constant $C_{SG}$ in inequality~\eqref{eq-poi}. This rate is independent of the initial condition.  For  reaction-diffusion equations, Theorem~\ref{thm-reaction} and Theorem~\ref{thm-2} prove that the optimal rate of convergence strongly depends on the initial conditions, which is natural for a chemical reaction, the reaction will converges quickly if the species are more concentrated at the beginning.  
%
\end{remark}
\begin{remark}
The result obtained in Theorem~\ref{thm-2} can be of course generalizes in the case when the generator $L$ is just an operator satisfying a spectral gap inequality on the domain considered $\Omega$. For example one can consider the case of a fractional Laplacian, a $p$-Laplacian,..., and others. (A problem may remain to prove the  existence of a non-negative solution in some of these cases.) 
\end{remark}

\begin{remark}
We recall also that while one can find in the literature some results on long time behaviour of solutions,  (as for example nice bounds in \cite{smoller}), generally it is a challenge to obtain the optimal bounds. Contrary to Theorem~\ref{thm-2}, the rate obtained in~\cite{dfm} for a general diffusion always depends on the spectral gap constant.  
\end{remark}
 \section{Study  of the general case} 
 \label{sec-5}
Let us consider now the general case for $q\geq 1$,
$$
\sum_{i=1}^q \alpha_i\A_i \rightleftharpoons \sum_{i=1}^q \beta_i\A_i.
$$
Assume now that $\Omega\in\R^n$ is  bounded domain.


Let consider non-negative initial  conditions $v_i^0\geq$ on $\Omega$, for all $1\leq i\leq q$. 
{\it A weak solution} of~\eqref{eq-def2}  on the time interval $I$ is $q$ measurable functions $\PAR{v_i}_{1\leq i\leq q}$ such that for all $t\in I$, $v_i(\cdot,t)\in L_1\PAR{\Omega}$, $G(v_1(\cdot, t),\cdots,v_q(\cdot, t))\in L_1\PAR{\Omega}$,  $$\int_0^t \NRM{G(v_1(\cdot, s),\cdots,v_q(\cdot, s))}_{L_1}ds<+\infty$$ and moreover for all $1\leq i\leq q$, $x\in\dR^n$ and $t>0$,
\begin{equation}
\label{eq-100}
   \quad v_i(t,x)=\Pt\PAR{v_i^0}+k_i(\beta_i-\alpha_i)\int_0^t\PT{t-s}\PAR{G(v_1(\cdot, s),\cdots,v_q(\cdot, s))}ds,
\end{equation}
which satisfies also for all $x\in\Omega$, $v_i(0,x)=v_i^0(x)$.

The result given here is a direct application of Rothe \cite[Theorem~4]{rothe}. The main interest of this result is to see why the solution remains non-negative and is defined on $[0,\infty)$ which generally can be quite surprising for a fully nonlinear parabolic equation; (see also \cite{smoller} for arguments based on a comparison principle). 

\begin{eprop}
\label{thm-eu}
Assume that there exist $1\leq i_0,j_0\leq q$ such that $\beta_{i_0}-\alpha_{i_0}>0$ and $\beta_{j_0}-\alpha_{j_0}<0$.

Then, for any non-negative bounded and measurable  initial condition $\PAR{v_i^0}_{1\leq i\leq q}$, there exists a non-negative weak solution of the system~$(\ref{eq-def2})$.
\end{eprop}
\begin{eproof}
We will give here just a sketch of the proof and refer to  \cite{rothe} and references therein to get more informations.   

\medskip

Let us see how solutions are bounded and non-negative. The idea is to solve a different problem : As for equation~(\ref{eq-ct}), for all $1\leq i\leq q$ we note by $C_{i,1}(t,x)=\Pt\PAR{v_i^0-\frac{k_i(\beta_i-\alpha_i)}{k_1(\beta_1-\alpha_1)}v_1^0}(x)$ and consider the  PDE
\begin{equation}
\label{eq-apt3}
\partial_t v_1=  L v_1 +\bar{F}(x,t,v_1),
\end{equation}
where 
\begin{multline*}
  \bar{F}(x,t,v_1)=k_1(\beta_1-\alpha_1)\\
  \PAR{\prod_{i=1}^q \PAR{\frac{k_i(\beta_i-\alpha_i)}{k_1(\beta_1-\alpha_1)}v_1+C_{i,1}(t,x)}_+^{\alpha_i} -\prod_{i=1}^q \PAR{\frac{k_i(\beta_i-\alpha_i)}{k_1(\beta_1-\alpha_1)}v_1+C_{i,1}(t,x)}_+^{\beta_i}},
\end{multline*}
and for $x\in \dR$, $(x)_+=\max\BRA{x,0}$.
%
%
%
%

By~\cite{rothe},~\eqref{eq-apt3} has an optimal and weak solution. Let us see why the solution is bounded and non-negative. Let $\bar{v}_1^0$ be a non-negative bounded initial condition and let  $\bar{v}_1$ a weak  solution of~(\ref{eq-apt3}) (with the same definition as for equation~(\ref{eq-def2}), replacing $G$ by $F$ in (\ref{eq-100})). 
 


Let for all $1\leq i\leq q$, 
\begin{equation}
\label{eq-dee}
\bar{v}_i(t,x)=\frac{k_i(\beta_i-\alpha_i)}{k_1(\beta_1-\alpha_1)}\bar{v}_1(t,x)+C_{i,1}(t,x).
\end{equation}
Then $(\bar{v}_i)_{1\leq i\leq q}$ is a solution of 
\begin{equation}
\label{eq-apt4}
\partial_t \bar{v}_i=  L \bar{v}_i +k_i(\beta_i-\alpha_i)\bar{G}(\bar{v}_1,\cdots,\bar{v}_q),
\end{equation}
where 
$$
\bar{G}(\bar{v}_1,\cdots,\bar{v}_q)={\prod_{i=1}^q \PAR{\bar{v}_i}_+^{\alpha_i} -\prod_{i=1}^q \PAR{\bar{v}_i}_+^{\beta_i}}.
$$
Let us multiply \eqref{eq-apt4} by $-(\bar{v_i})_-:=\min\BRA{\bar{v}_i,0}$. After integration, we obtain
\begin{equation}
\label{eq-pos}
  \quad\frac{d}{dt}\int \frac 1 2 ((\bar{v}_i)_-)^2d\mu=\int (\bar{v}_i)_-L\PAR{(\bar{v}_i)_ -}d\mu -k_i(\beta_i-\alpha_i)\int (\bar{v}_i)_-\bar{G}(\bar{v}_1,\cdots,\bar{v}_q)d\mu,
\end{equation}
where $\int (\bar{v}_i)_-L\PAR{(\bar{v}_i)_ -}d\mu=\int \bar{v}_iL\PAR{(\bar{v}_i)_ -}d\mu\leq0$. 

On the set $\BRA{\bar{v}_i\leq0}$, we have 
$$
k_i(\beta_i-\alpha_i)\bar{G}(\bar{v}_1,\cdots,\bar{v}_q)=k_i(\beta_i-\alpha_i)\bar{G}(\bar{v}_1,\cdots,0,\cdots,\bar{v}_q),
$$
 where we put $0$ at the position $i$. Since for all $j$, $(\bar{v}_j)_+\geq0$ then it is not difficult to see that in all cases  $k_i(\beta_i-\alpha_i)\bar{G}(\bar{v}_1,\cdots,0,\cdots,\bar{v}_q)\geq0$. Which gives that 
$$
\frac{d}{dt}\int \frac 1 2 ((\bar{v}_i)_-)^2d\mu\leq0.
$$ 
Since at time $t=0$, $\int ((\bar{v}_i^0)_-)^2d\mu =0$ then for all $t\geq 0$, $\bar{v}_i(t)\geq0$ almost everywhere. 

\medskip

Assume that $\beta_1-\alpha_1>0$. Then, since  the  solutions are non-negative, we get the following global estimate of the solution,  
\begin{multline*}
  0\leq \bar{v}_1 =\frac{k_1(\beta_1-\alpha_1)}{k_{j_0}(\beta_{j_0}-\alpha_{j_0})}\bar{v}_{j_0}+\\
\Pt\PAR{\bar{v}_1^0-\frac{k_1(\beta_1-\alpha_1)}{k_{j_0}(\beta_{j_0}-\alpha_{j_0})}\bar{v}_{j_0}^0}\leq \NRM{\bar{v}_1^0}_\infty+\ABS{\frac{k_1(\beta_1-\alpha_1)}{k_{j_0}(\beta_{j_0}-\alpha_{j_0})}}\NRM{\bar{v}_{j_0}^0}_\infty. 
\end{multline*}
If $\beta_1-\alpha_1>0$ does not hold, we use $i_0$ instead $j_0$ to get the same result with $j_0$. The last estimate proves that  $ \bar{v}_1$ is bounded and then the solution is defined on $[0,\infty)$.


The same method as above in~\eqref{eq-pos} proves that for all $0\leq i\leq q$, $\bar{v}_i\geq 0$. This implies that $\bar{G}\PAR{\bar{v}_1,\cdots,\bar{v}_q}={G}\PAR{\bar{v}_1,\cdots,\bar{v}_q}$ and then $(\bar{v}_i)_{1\leq i\leq q}$ is also a non-negative weak solution of~\eqref{eq-def2} which finished the proof of the existence. 
\end{eproof}

\begin{remark}
This restriction on parameters $(\alpha_i,\beta_i)$ is natural in the context of a chemical  reaction by the the principle of conservation of mass by Lavoisier. 
\end{remark}

The following corollary is a direct consequence of~\eqref{eq-ct}.

\begin{ecor}
\label{cor-ine}
Assume that there exist $1\leq i_0,j_0\leq q$ such that $\beta_{i_0}-\alpha_{i_0}>0$ and $\beta_{j_0}-\alpha_{j_0}<0$.

Let $(v_i)_{1\leq i\leq q}$ be a solution of~\eqref{eq-def2}.  Then for all $1\leq i\leq q$ such that $\beta_{i}-\alpha_{i}>0$ one gets for all $t\geq0 $ and $x\in\Omega$ $:$
$$
0\leq v_i(t,x)\leq k_i(\beta_{i}-\alpha_{i})\min\BRA{\Pt\PAR{\frac{v_i^0}{k_i(\beta_{i}-\alpha_{i})}-\frac{v_j^0}{k_j(\beta_{j}-\alpha_{j})}}(x),\,\,\beta_{j}-\alpha_{j}<0}
$$
and for all $1\leq i\leq q$ such that $\beta_{i}-\alpha_{i}<0$ $:$
$$
0\leq v_i(t,x)\leq k_i\ABS{\beta_{i}-\alpha_{i}}\min\BRA{\Pt\PAR{\frac{v_i^0}{k_j(\beta_{j}-\alpha_{j})}-\frac{v_i^0}{k_i(\beta_{i}-\alpha_{i})}}(x),\,\,\beta_{j}-\alpha_{j}>0}.
$$
In particular when initial conditions are bounded, solutions of~$\eqref{eq-def}$ are also bounded with an explicit upper bound. 
\end{ecor}

\begin{ethm}
\label{eq-gene}
Assume that the semigroup $(\Pt)_{t\geq 0}$  satisfies a spectral gap inequality~\eqref{eq-sg} with respect to the invariant probability measure $\mu$. Assume also that for all $1\leq i \leq q$, $\alpha_i\beta_i=0$.  

Let $\PAR{v_i^0}_{1\leq i\leq q}$ be a  non-negative bounded initial condition. We assume furthermore  that for all $1\leq i\leq q$, $\int v_i^0 d\mu>0$. 

Let   $(s_i)_{1\leq i\leq q}$ be the steady state given by Lemma~\ref{lem-sta}. Then for all $1\leq i\leq q$, one gets 
\begin{equation}
\label{eq-fin}
\sqrt{\int \PAR{v_i-s_i}^2 d\mu}\leq K\exp\PAR{-\min\BRA{a,M}t},
\end{equation}
where $a>0$ depends on $\alpha_i$, $\beta_i$ and $C_{SG}$, and  $M,K>0$  depend on the initial conditions.
\end{ethm}

\begin{eproof}
The idea is almost the same as in Theorem~\ref{thm-2} except that we do not obtain the optimal rate. 

Assume that $\beta_1-\alpha_1>0$, the opposite case could be treated in the similar way.  Equation~\eqref{eq-apt} applied for $j=1$ reads
\begin{equation*}
\partial_t v_1=L v_1 +F(t,x,v_1),
\end{equation*}
where 
\begin{multline*}
   F(t,x,y)=k_1(\beta_1-\alpha_1)\times\\
\PAR{\prod_{i=1}^q \PAR{\frac{k_i(\beta_i-\alpha_i)}{k_1(\beta_1-\alpha_1)}y+C_{i,1}(t,x)}^{\alpha_i} -\prod_{i=1}^q \PAR{\frac{k_i(\beta_i-\alpha_i)}{k_1(\beta_1-\alpha_1)}y+C_{i,1}(t,x)}^{\beta_i}}
\end{multline*}
 and functions $C_{i,1}(t,x)$ are defined in~\eqref{eq-ct}. By the ergodicity properties of the semi-group, equation~\eqref{eq-ergodi}, implies that in $L^2(d\mu)$, 
$$
\lim_{t\rightarrow \infty}C_{i,1}(t,x)=\int \PAR{v_i^0-\frac{k_i(\beta_i-\alpha_i)}{k_1(\beta_1-\alpha_1)}v_1^0}d\mu:=C_{i,1}^\infty.
$$
Denote by $F_\infty(y)$ the limit of $F(t,x,y)$ when $t$ goes to infinity. Note that $F_\infty$ does not depend on $x\in\Omega$. 
Then, one gets since $\int v_1Lv_1d\mu \leq0$, 
\begin{eqnarray} 
\frac{d}{dt}\frac{1}{2}\int\PAR{v_1-s_1}^2 d\mu\leq 
k_1(\beta_1-\alpha_1)\int \PAR{v_1-s_1}(F(t,\cdot,v_1)-F_\infty(v_1))d\mu\nonumber\\\quad\quad\quad\quad\quad\quad\quad\quad\quad\quad\quad+k_1(\beta_1-\alpha_1)\int\PAR{v_1-s_1}F_\infty(v_1)d\mu.\label{eq-bof}
\end{eqnarray}
Note that $F_\infty$ is equal to the polynomial function  $F$ defined in~\eqref{eq-deff} where in the definition of $C_{j,i}$ in~\eqref{eq-defcji},  $v_i(0)$ is replaced by the mean value of initial conditions $\int v_i^0d\mu$. Let us set 
$$
M_1(t,x)=k_1(\beta_{1}-\alpha_{1})\min\BRA{\Pt\PAR{\frac{v_1^0}{k_1(\beta_{1}-\alpha_{1})}-\frac{v_j^0}{k_j(\beta_{j}-\alpha_{j})}}(x),\,\,\beta_{j}-\alpha_{j}<0},
$$
and for its limit as $t$ goes to $\infty$ 
$$
M_1^\infty=k_1(\beta_{1}-\alpha_{1})\min\BRA{\int\PAR{\frac{v_1^0}{k_1(\beta_{1}-\alpha_{1})}-\frac{v_j^0}{k_j(\beta_{j}-\alpha_{j})}}d\mu,\,\,\beta_{j}-\alpha_{j}<0}.
$$

As it was shown in the proof in Theorem~\ref{thm-reaction},  $F_\infty(X)=(X-s_1)Q(X)$ with $Q(X)<0$ for all $X\in(0,M_1^\infty)$. Now since for all $1\leq i\leq q$, $\alpha_i\beta_i=0$ then $s_1$ is a simple root of the polynomial function $F_\infty$. It implies that  that  $Q(0)<0$ and  $Q(M_1^\infty)<0$. Then by continuity of $Q$ there exist $\epsilon,\eta>0$ such that $Q(X)\leq-\epsilon$ for all $X\in[0,M_1^\infty+\eta]$.

For the second term in \eqref{eq-bof}, we get 
\begin{eqnarray*}
   k_1(\beta_1-\alpha_1)\int\PAR{v_1-s_1}F_\infty(v_1)d\mu\leq
-k_1(\beta_1-\alpha_1)\epsilon\int_\BRA{v_1\leq M_1^\infty+\eta}\PAR{v_1-s_1}^2d\mu\\\quad\quad\quad\quad\quad\quad\quad\quad\quad\quad+k_1(\beta_1-\alpha_1)\int_\BRA{v_1> M_1^\infty+\eta}\PAR{v_1-s_1}F_\infty(v_1)d\mu,
\end{eqnarray*}
and then for some constant $K$ depending on initial conditions
\begin{multline*}
k_1(\beta_1-\alpha_1)\int\PAR{v_1-s_1}F_\infty(v_1)d\mu\leq\\-k_1(\beta_1-\alpha_1)\epsilon\int\PAR{v_1-s_1}^2d\mu+K\mu\BRA{v_1> M_1^\infty+\eta}.
\end{multline*}
Corollary~\ref{cor-ine} implies that $M_1(t,\cdot)\geq v_1$ and then Markov inequality gives 
\begin{multline}
\label{eq-ine54}
  k_1(\beta_1-\alpha_1)\int\PAR{v_1-s_1}F_\infty(v_1)d\mu\leq\\-k_1(\beta_1-\alpha_1)\epsilon\int\PAR{v_1-s_1}^2d\mu+K\var{\mu}{M_1(t,\cdot)}.
\end{multline}
Since for $q$ measurable  functions $g_i\in L_2(\mu)$ one has
$$
\var{\mu}{\min\BRA{g_i,\,\,1\leq i\leq q}}\leq \frac{1}{2}\sum_{i=1}^q\var{\mu}{g_i},
$$
so the last term of~\eqref{eq-ine54} gives 
\begin{equation}
\label{eq-finpeut}
   \quad k_1(\beta_1-\alpha_1)\int\PAR{v_1-s_1}F_\infty(v_1)d\mu\leq-k_1(\beta_1-\alpha_1)\epsilon\int\PAR{v_1-s_1}^2d\mu+K' e^{-\frac{1}{C_{}SG}t},
\end{equation}
where $K'$ is an another constant depending on initial conditions.

Let us note $F(t,x,y)=\sum_{i=1}^\gamma K_{t,i,x}y^i$ and $F_\infty(y)=\sum_{i=1}^\gamma K_{\infty,i}y^i$ where we note $\gamma=\max\BRA{\sum_{i=1}^q \beta_i,\sum_{i=1}^q \alpha_i}$. The first term of~\eqref{eq-bof} gives 
\begin{multline*}
   \int \PAR{v_1-s_1}(F(t,\cdot,v_1)-F_\infty(v_1))d\mu\\\leq\sqrt{\int \PAR{v_1-s_1}^2d\mu}\sum_{i=1}^\gamma\sqrt{\int v_1^{2i}(K_{t,i,x}-K_{\infty,i})^2d\mu}. 
\end{multline*}
Let us consider one of them, with $1\leq i\leq q$, one has 
$$
\int v_1^{2i}(K_{t,i,x}-K_{\infty,i})^2d\mu\leq\sqrt{\int v_1^{4i}d\mu}\sqrt{\int (K_{t,i,x}-K_{\infty,i})^4d\mu}.
$$
There exist some sets $\Gamma_{i,j}$ and $\Delta_i$ and constants $\mu_{i,j,k}\in\dR$, $\gamma_{i,j,k}\in\dN$ such that 
$$
K_{t,i,x}=\sum_{j\in\Delta_i}\prod_{k\in\Gamma_{i,j}}\mu_{i,j,k}\PAR{C_{i,1}(t,x)}^{\gamma_{i,j,k}}.
$$
Then for some constant $K_i>0$, $p_{i,j,k}\geq 2$ and $q_{i,j,k}>0$,
$$
\int (K_{t,i,x}-K_{\infty,i})^4d\mu\leq K_i \sum_{j\in\Delta_i}\prod_{k\in\Gamma_{i,j}}\PAR{\int (\PAR{C_{i,1}(t,\cdot)}^{\gamma_{i,j,k}}-\PAR{C_{i,1}^\infty}^{\gamma_{i,j,k}})^{p_{i,j,k}}d\mu}^{q_{i,j,k}}.
$$
Since the initial conditions are bounded and $p_{i,j,k}\geq 2$, one gets for some another constant $K$  
$$
\int (\PAR{C_{i,1}(t,\cdot)}^{\gamma_{i,j,k}}-\PAR{C_{i,1}^\infty}^{\gamma_{i,j,k}})^{p_{i,j,k}}d\mu\leq K\int \PAR{{C_{i,1}(t,\cdot)}-{C_{i,1}^\infty}}^2d\mu,
$$
and then spectral gap inequality gives for some  $K'$, 
$$
\int (\PAR{C_{i,1}(t,\cdot)}^{\gamma_{i,j,k}}-\PAR{C_{i,1}^\infty}^{\gamma_{i,j,k}})^{p_{i,j,k}}d\mu\leq K'e^{-\frac{1}{C_{SG}}t}.
$$
Thus we have proved that there exits  $\gamma_i>0$ and $R_i>0$ depending on initial conditions and $C_{SG}$ such that 
$$
\int (K_{t,i,x}-K_{\infty,i})^4d\mu \leq R_i e^{-\gamma_i t}. 
$$
All of these estimates give for some $\alpha>0$ depending on $\alpha_i$, $\beta_i$ and $C_{SG}$ and  $R>0$ depending on initial conditions the following bound
$$
k_1(\beta_1-\alpha_1)\int \PAR{v_1-s_1}(F_t(v_1)-F_\infty(v_1))d\mu\leq\sqrt{\int \PAR{v_1-s_1}^2d\mu}  \,R e^{-\alpha t}. 
$$
If $\phi(t)=\sqrt{\int \PAR{a_1-s_1}^2d\mu}$, the equation~\eqref{eq-bof} becomes, by the previous equation and~\eqref{eq-finpeut},
$$
\phi'(t)\leq -k_1(\beta_1-\alpha_1)\epsilon\phi(t)+R e^{-\alpha t}+K' e^{-\frac{1}{C_{}SG}t},
$$
which finishes the proof. 
\end{eproof}

\begin{remark}
\begin{itemize}
\item 
One can generalize the last theorem in the following way without assuming that for all $i$, $\alpha_i \beta_i=0$. 
Let us consider $i^+$ and $j^-$ such that 
$$
\sup_{j,\,s.t.\,\beta_j-\alpha_j>0} \BRA{-\frac{\int v_j^0d\mu}{\beta_j-\alpha_j}}=-\frac{\int v_{i^+}^0d\mu}{\beta_{i^+}-\alpha_{i^+}}
$$
and 
$$
\sup_{j,\,s.t.\,\beta_j-\alpha_j<0} \BRA{-\frac{\int v_j^0d\mu}{\beta_j-\alpha_j}}=-\frac{\int v_{j^-}^0d\mu}{\beta_{j^-}-\alpha_{j^-}}. 
$$
Assume only that $\alpha_{i^+} \beta_{i^+}=\alpha_{j^-} \beta_{j^-}=0$. Then the computation for the species $i^+$ and $j^-$ are the same as in the proof of Theorem~\ref{eq-gene}.
\item On can also generalize in assuming that the initial conditions $(v_i^0)_{1\leq i\leq q}$ are in $L^q(d\mu)$ for some $q>1$ instead of  $L^{\infty}(d\mu)$. In that case the proof will be more technical.  
\end{itemize}
\end{remark}

\subsubsection*{Acknowledgment}
\noindent We would like to thank the anonymous referee for her/his valuable comments on the previous version. \\
This research was supported in part by the ANR project EVOL and by EPSRC grant EP/D05379X/1.  The first author thanks the members of UMPA at the \'Ecole Normale Sup\'erieure de Lyon for their kind hospitality.


\medskip\noindent

\noindent
Ceremade\\
Universit\'e Paris-Dauphine\\
Place du Mar\'echal De Lattre De Tassigny\\
75116 Paris - France\\
gentil@ceremade.dauphine.fr\\

\noindent
Institut de Math\'ematiques de Toulouse\\
Universit\'e de Toulouse\\
31062 Toulouse - France\\
and \\
Imperial College London\\
South Kensington Campus\\
London SW7 2AZ, UK\\
Email: b.zegarlinski@imperial.ac.uk \\

\medskip

\noindent

\noindent
\end{document}